\documentclass[12pt]{amsart}

\usepackage{amscd}
\usepackage{amssymb,amsmath,amsthm}
\theoremstyle{plain}
\newtheorem{thm}{Theorem}[section]
\newtheorem{lem}[thm]{Lemma}

\newtheorem{defn-lem}[thm]{Definition-Lemma}

\theoremstyle{definition}
\newtheorem{defn}[thm]{Definition}

\newtheorem{rem}[thm]{Remark}

\def\md #1#2#3#4#5 {\left(
                        \begin{matrix}
             #1 & #2 \\
             #3 & #4
                        \end{matrix}
                      \right)- #5}

\def\ma #1#2#3#4 {\left(
                        \begin{matrix}
             #1 & #2 \\
             #3 & #4
                        \end{matrix}
                      \right)}
\def\mu (#1) {\mathcal{M}#1}

\def\Ad{\operatorname{Ad}}
\def\id{\operatorname{id}}
\def\Aut{\operatorname{Aut}}
\def\KK{\operatorname{KK}}
\newcommand{\mc}{\mathcal}
\begin{document}
\title [Proper asymptotic unitary equiv. and projection lifting ]
       {Proper asymptotic unitary equivalence in $\KK$-theory and projection lifting from the corona algebra }

\begin{abstract}
In this paper we generalize the notion of essential codimension of Brown, Douglas, and Fillmore using $\KK$-theory and prove a result which asserts that there is a unitary of the form `identity + compact' which gives the unitary equivalence of two projections if the `essential codimension' of two projections vanishes for certain $C\sp*$-algebras employing the proper asymptotic unitary equivalence of $\KK$-theory found by M. Dadarlat and S. Eilers. We also apply our result to study the projections in the corona algebra of $C(X)\otimes B$ where $X$ is $[0,1]$, $(-\infty, \infty)$, $[0,\infty)$, and $[0,1]/\{0,1\}$.
\end{abstract}

\author { Hyun Ho \quad Lee }

\address {Department of Mathematics\\
         Seoul National University\\
         Seoul, South Korea 151-747 }
\email{hadamard@snu.ac.kr}


\date{}
\maketitle

\section{Introduction}
When two projections $p$ and $q$ in $B(H)$, whose difference is compact, are given, an integer $[p:q]$ is defined as the Fredholm index of $v^*w$ where $v,w$ are isometries on $H$ with $vv^*=p$ and $ww^*=q$. This number is called the essential codimension because it gives the codimension of $p$ in $q$ if $p\leq q$ \cite{BDF}. A modern interpretation of this essential codimension is provided using the Kasparov group $\KK(\mathbb{C},\mathbb{C})$. Indeed, a $*$-homomorphism from $\mathbb{C}$ to $B(H)$ is determined by the image of $1$ which is a projection. Thus we can associate to the essential codimension  a Cuntz pair.
An important result of the essential codimension is the following: $[p:q]=0$ if and only if there is a unitary $u$ of the form identity + compact such that $upu^*=q$. Motivated by this result, Dadarlat and Eilers defined a new equivalence relation on $\KK$-group \cite{DE}. When  $\pi,\sigma: A \to \mathcal{L}(E)$ are two representations, with $E$ is a Hilbert $B$-module, we say $\pi$ and $\sigma$ are properly asymptotically unitarily equivalent and write $\pi \approxeq \sigma$ if there is a continuous path of unitaries $u:[0,\infty) \to \mathcal{U}(\mathcal{K}(E)+\mathbb{C}1_E)$, $u=(u_t)_{t\in[0,\infty)}$, such that
\begin{itemize}
\item $lim_{t\to \infty}\| \sigma(a)-u_t \pi(a)u_t^*\|=0$ for all $a \in A$,
\item $ \sigma(a)-u_t \pi(a)u_t^* \in \mathcal{K}(E)$ for all $t\in [0,\infty )$, and $a\in A$.
\end{itemize}

Note that the word `proper' reflects the fact that implementing unitaries are of the form `identity+compact'. The main result of them is \cite[Theorem 3.8]{DE} which asserts that if $\phi, \psi: A \to M(B\otimes K(H))$ is a Cuntz pair of representations, then the class $[\phi,\psi]$ vanishes in $\KK(A,B)$ if and only if there is another representation $\gamma:A \to M(B\otimes K(H))$ such that $\phi\oplus\gamma \approxeq \psi\oplus \gamma$. When $B=\mathbb{C}$, which corresponds to $K$-homology, the result is improved as a non-stable version. In fact, if $(\phi,\psi)$ is a Cuntz pair of faithful, non-degenerate representations from $A$ to $B(H)$ such that both images do not contain any non-trivial compact operator, then the cycle $[\phi,\psi]=0 $ in $\KK(A,\mathbb{C})$ if and only if $\phi \approxeq \psi$ \cite[Theorem 3.12]{DE}. This fits nicely with the above aspect of the essential codimension. An abstract version of this is proved given a Cuntz pair of absorbing representations (See Theorem \ref{T:DE}). Thus the proper asymptotic unitary equivalence must be the right notion and tool for further developments of the non-stable K-theory. Our intrinsic interest lies in when this non-stable version of proper asymptotic unitary equivalence happens as shown in $K$-homology case. We show a similar result for $K$-theory. In fact, we prove that if $(\phi,\psi) $ is a Cuntz pair of faithful representations from $\mathbb{C} \to M(B\otimes K)$ whose images are not in $B\otimes K$, then $[\phi,\psi]=0$ in $K(B)$ if and only if $\phi \approxeq \psi$ provided that $B$ is non-unital, separable, purely infinite simple $C\sp*$-algebra such that $M(B)$ has real rank zero (See Theorem \ref{T:BDF}).

Besides our intrinsic interest, Theorem \ref{T:BDF} was motivated by the projection lifting problem from the corona algebra to the multiplier algebra of a $C\sp*$-algebra of the form $C(X)\otimes B$. To lift a projection from a quotient algebra to a projection  has been a fundamental question related to K-theory (See \cite{Da}). We show that  a projection in the corona algebra is `locally' liftable to a projection in the multilpler algebra  but not `globally' in general. In other words, it can be represented by finitely many  projection valued functions so that  its discontinuities are described in terms of  Cuntz pairs. They give rise to $K$-theoretical obstructions. We show that these discontinuities  can be resolved if corresponding $K$-theoretical terms are vanishing. In this process, the crucial point of proper asymptotic unitary equivalence is exploited as a key step (See Theorem \ref{T:lifting}).

\section{Proper asymptotic unitary equivalence}
Let $E$ be a (right) Hilbert $B$-module. We denote $\mathcal{L}(E,F)$ by the $C\sp{*}$-algebra of adjointable, bounded  operators from $E$ to $F$. The ideal of `compact' operators from $E$ to $F$ is denoted by $\mathcal{K}(E,F)$. When $E=F$, we write $\mathcal{L}(E)$ and $\mathcal{K}(E)$ instead of  $\mathcal{L}(E,E)$ and $\mathcal{K}(E,E)$. Throughout the paper, $A$ is a separable $C\sp{*}$-algebra, and all Hilbert modules are assumed to be countably generated over a separable $C\sp*$-algebra. We use the term representation for a $*$-homomorphism from $A$ to $\mathcal{L}(E)$. We let $H_B$ be the standard Hilbert module over $B$ which is $H\otimes B $ where $H$ is a separable infinite dimensional Hilbert space. We denote $M(B)$ by the multiplier algebra of $B$. It is well-known that $\mathcal{L}(H_B)=M(B\otimes K)$ and $\mathcal{K}(H_B)=B\otimes K$ where $K$ is the $C\sp{*}$-algebra of the compact operators on $H$ \cite{Kas80}.
\begin{defn}\cite[Definition 2.1]{DE}
Let $\pi,\sigma$ be two representations from $A$ to $E$ and $F$ respectively. We say $\pi$ and $\sigma$ are approximately unitarily equivalent and write $\pi \sim \sigma$, if there exists a sequence of unitaries $u_n \in \mathcal{L}(E,F)$ such that for any $a\in A$
\begin{itemize}
\item[(i)] $lim_{n\to \infty}\| \sigma(a)-u_n \pi(a)u_n^*\|=0$,\\
\item[(ii)] $ \sigma(a)-u_n \pi(a)u_n^* \in \mathcal{K}(F)$ for all $n$.
\end{itemize}
\end{defn}
\begin{defn}\cite[Definition 2.5]{DE}
A representation $\pi:A \to \mathcal{L}(E)$ is called absorbing if $\pi \oplus \sigma \sim \pi$ for any representation $\sigma :A \to \mathcal{L}(F)$.
\end{defn}
  We say that $\pi$ and $\sigma$ are asymptotically unitarily equivalent, and write $\pi \underset{\text{asym}}{\sim}\sigma$ if there is a unitary valued norm continuous map $u:\left[0,\infty \right) \to \mathcal{L}(E,F)$ such that $t \to \sigma(a)- u_t\pi(a)u_t^* $ lies in $C_0([0,\infty))\otimes \mathcal{K}(E)$ for any $a \in A$, or if
\begin{itemize}
\item[(i)]$lim_{t\to \infty}\| \sigma(a)-u_t \pi(a)u_t^*\|=0$,\\
\item[(ii)] $ \sigma(a)-u_t \pi(a)u_t^* \in \mathcal{K}(F)$ for all $t\in [0,\infty)$.
\end{itemize}
If $\pi:A \to \mathcal{L}(E)$ is a representation, we define $\pi^{(\infty)}:A \to \mathcal{L}(E^{(\infty)})$ by as $\pi^{(\infty)}=\pi\oplus \pi\oplus \cdots$ where $E^{(\infty)}=E\oplus E\oplus \cdots$.
\begin{lem}\label{L:absorbing}
 Let $\psi$ be an absorbing representation, and $\phi$ be a representation, of a separable $C\sp{*}$-algebra $A$ on the standard Hilbert $B$-module $H_B$. Then there exists a sequence of isometries $\{v_n\} \subset \mathcal{L}(H_B^{(\infty)},H_B)$ such that for each $a \in A$
\begin{align*}
   &v_n \phi^{(\infty)}(a) -\psi(a) v_n \in \mathcal{K}(H_B^{(\infty)}, H_B), \\
   &\| v_n \phi^{(\infty)}(a) -\psi(a) v_n  \| \to 0 \quad \text{as $n \to \infty$},\\
   &v_j^* v_i=0 \quad \text{for $i \neq j$}.
\end{align*}
\end{lem}
\begin{proof}

Let $S_i, \, i=1,2,3 \cdots $, be a sequence of isometries of $\mc{L}(H_B)$ such that $S_i^*S_j=0, \,i\neq j$, and $\sum_iS_iS_i^*=1$ in the strict topology. Let $\phi_{\infty}(a)=\sum_i S_i \phi(a) S_i^*$. Since $\psi$ is absorbing, there is a unitary $U \in \mathcal{L}(H_B, H_B)$ such that
\begin{equation}\label{E:(1)}
   U^* \psi(a)U- \phi_{\infty}(a)  \in \mc{K}(H_B) \quad a\in A.
\end{equation}
Define $T:H_B^{(\infty)} \to H_B$ by $T=(S_1,S_2,\cdots )$. Then  $$\phi_{\infty}(a)=T \phi^{(\infty)}(a)T^*. $$
Thus equation (\ref{E:(1)}) is rewritten as
\begin{equation}\label{E:(2)}
  T^*U^* \psi(a)UT- \phi^{(\infty)}(a)  \in \mc{K}(H_B^{(\infty)}) \quad a\in A.
 \end{equation}
 If we identify $\phi^{(\infty)}$ as $(\phi^{(\infty)})^{(\infty)}$, there is a partition $N_i$, $i=1,2,3,\cdots$, of $\mathbb{N}$ so that we generate a sequence of isometries $v_i \in \mc{L}(H_B^{(\infty)},H_B)$ from $UT=(US_1,US_2,\cdots,)$. More concretely, if we let $\nu_i: N_i \to \mathbb{N}$ be bijections, we can define
 $v_i=(US_{\nu_i^{-1}(1)},US_{\nu_i^{-1}(2)},\cdots)$. It is easily checked that $v_iv_j^*=0$ for $i\neq j$.
 Equation (\ref{E:(2)}) implies that
 \begin{align*}
 & v_i^*\psi(a) v_i -\phi^{(\infty)}(a) \in \mc{K}(H_B^{(\infty)}), \\
 &\left\| v_i^*\psi(a) v_i -\phi^{(\infty)}(a) \right\| \to 0 \quad \text{as}\, i \to \infty.
 \end{align*}
 Finally, our claim follows from
 \[
 \begin{split}
 &(v_n \phi^{(\infty)}(a) -\psi(a) v_n)^*(v_n \phi^{(\infty)}(a) -\psi(a) v_n)=\phi^{(\infty)}(a^*)(\phi^{(\infty)}(a)-v_n^*\psi(a)v_n)\\
 &+(\phi^{(\infty)}(a^*)-v_n^*\psi(a) v_n)\phi^{(\infty)}(a)-(\phi^{(\infty)}(a^*a)-v_n^*\psi(a^*a)v_n).
 \end{split}
 \]
\end{proof}

\begin{lem}\cite[Lemma 2.6]{DE}\label{L:asym}
Let $\pi:A \to \mathcal{L}(E)$ and $\sigma:A \to \mathcal{L}(F)$ be two representations. Suppose that there is a sequence of
isometries $v_i:F^{(\infty)} \to E$ such that for  $a\in A$
$$ v_i\sigma ^{(\infty)}(a)-\pi(a)v_i \in \mathcal{K}(F^{(\infty)}, E), \quad \lim_{i \to \infty}\|v_i\sigma ^{(\infty)}(a)-\pi(a)v_i\|\to 0, $$
and $v_j^*v_i=0$ for $i\neq j$. Then $\pi \oplus\sigma \underset{\text{asym}}{\sim}\pi$.
\end{lem}
 We say $\phi:A \to B(H)$ is admissible if $\phi$ is faithful, non-degenerate, and $\phi(A)\cap K = \{0 \}$. The main result in \cite {Voi}
states that any pair of admissible representations $\phi$ and $\psi$ satisfies that $\phi \sim \psi$. Dadarlat and Eilers proved a much stronger version which states that any pair of admissible representations $\phi$ and $\psi$ satisfies  $\phi \underset{\text{asym}}{\sim} \psi$ \cite[Theorem 3.11]{DE}. Since the admissible representation is absorbing, the following result is the appropriate generalization of Voiculescu's result.
\begin{thm}\label{T:Voiculescu}
If two representations $\psi$, $\phi$ of a separable $C\sp{*}$-algebra $A$ on the standard Hilbert $B$-module $H_B$ are absorbing,
then we have $\phi \underset{\text{asym}}{\sim} \psi$.
\end{thm}
\begin{proof}
By Lemma \ref{L:absorbing} and Lemma \ref{L:asym}, we have $\psi\oplus\phi  \underset{\text{asym}}{\sim} \psi$, and the proof is complete by symmetry.
\end{proof}

\begin{defn}
Let $\phi$ be a representation from $A$ to $M(B\otimes K)$. Then we define a $C\sp{*}$-algebra by
$$D_{\phi}(A,B)=\{x\in M(B\otimes K)\mid x\phi(a)-\phi(a)x \in B \otimes K, \quad a \in A \}.$$
\end{defn}

\begin{lem}
If $M(B\otimes K)$ has real rank zero, then $D_{\phi}(\mathbb{C},B)$ has real rank zero for any representation $\phi:\mathbb{C} \to M(B\otimes K)$.
\end{lem}
\begin{proof}
The proof of the lemma is essentially based on the argument due to Brown and Pedersen \cite{BP}.

Note that any representation $\phi: \mathbb{C} \to M(B\otimes K)$ is determined by $\phi(1)$, which is a projection in $M(B\otimes K)$. Say $\phi(1)=p$. Then we see that $D_{\phi}(\mathbb{C},B)=\{x\in M(B\otimes K)\mid x p- p x \in B \otimes K \}$.

To show $D_{\phi}(\mathbb{C},B)$ has real rank zero, it is enough to show any self-adjoint element in $D_{\phi}(\mathbb{C},B)$ is approximated by a self-adjoint, invertible element. Let $x$ be a self-adjoint element. Using the obvious matrix notation
\[x=\left(
                        \begin{matrix}
             a & c \\
             c^* & b
                        \end{matrix}
                      \right),\]
$x p- p x \in B\otimes K$ implies that $c$ is `compact', i.e., it is in $B\otimes K$.
Since $M(B\otimes K)$ has real rank zero, $pM(B\otimes K)p$ and $(1-p)M(B\otimes K)(1-p)$ have real rank zero. Given $\epsilon>0$ we can find $b_0$ invertible in $(1-p)M(B\otimes K)(1-p)$ with $b_0=b_0^*$ and $\|b-b_0 \|<\epsilon$. Then considering $a-cb_0^{-1}c^*$, we can find $a_0$ in $pM(B\otimes K)p$ with $a_0=a_0^*$ and $\| a-a_0\|<\epsilon$, such that $a_0-cb_0^{-1}c^*$ is invertible in $pM(B\otimes K)p$. Then
$\left(\begin{matrix}
             p & cb_0^{-1} \\
             0 & 1-p
                        \end{matrix}
                      \right) $, $\left(
                        \begin{matrix}
             p & 0 \\
             b_0^{-1}c^* & 1-p
                        \end{matrix}
                      \right)$ are in $D_{\phi}(\mathbb{C},B)$ since $ cb_0^{-1}$ is `compact'. Thus
\[x_0=\left(
                        \begin{matrix}
             a_0 & c \\
             c^* & b_0
                        \end{matrix}
                      \right)= \left(
                        \begin{matrix}
             p & cb_0^{-1} \\
             0 & 1-p
                        \end{matrix}
                      \right)\left(
                        \begin{matrix}
             a_0-cb_0^{-1}c^* & 0  \\
             0 & b_0
                        \end{matrix}
                      \right)\left(
                        \begin{matrix}
             p & 0 \\
             b_0^{-1}c^* & 1-p
                        \end{matrix}
                      \right)\]
is invertible in $D_{\phi}(\mathbb{C},B)$. Evidently $\|x-x_0 \| <\epsilon$, so we are done.
\end{proof}

Let us recall the definition of Kasparov group $\KK(A,B)$. We refer the reader to \cite{Kas81} for the general introduction of the subject. A $\KK$-cycle is a triple $(\phi_0,\phi_1,u)$, where $\phi_i:A \to \mathcal{L}(E_i)$ are representations and $u \in \mathcal{L}(E_0,E_1)$ satisfies that
\begin{itemize}
\item[(i)] $u\phi_0(a)-\phi_1(a)u \in \mathcal{K}(E_0,E_1)$,
\item[(ii)]$\phi_0(a)(u^*u-1)\in \mathcal{K}(E_0)$, $\phi_1(a)(uu^*-1)\in \mathcal{K}(E_1)$.
\end{itemize}
The set of all $KK$-cycles will be denoted by $\mathbb{E}(A,B)$.  A cycle is degenerate if
\[ u\phi_0(a)-\phi_1(a)u=0, \quad\phi_0(a)(u^*u-1)=0,\quad \phi_1(a)(uu^*-1)=0.\]
An operator homotopy through $KK$-cycles is a homotopy $(\phi_0,\phi_1,u_t)$, where the map $t \to u_t$ is norm continuous.
The equivalence relation $\underset{\text{oh}}{\sim}$ is generated by operator homotopy and addition of degenerate cycles up to unitary equivalence.
Then $\KK(A,B)$ is defined as the quotient of $\mathbb{E}(A,B)$ by $\underset{\text{oh}}{\sim}$. When we consider  non-trivially graded $C\sp{*}$-algebras, we define a triple $(E,\phi,F)$, where $\phi:A \to \mathcal{L}(E)$ is a graded representation, and $F\in \mathcal{L}(E)$ is of odd degree such that $F\phi(a)-\phi(a)F$, $(F^2-1)\phi(a)$, and $(F-F^*)\phi(a)$ are all in $\mathcal{K}(E)$ and call it a Kasparov $(A,B)$-module. Other definitions like degenerate cycle and operator homotopy are defined in  similar ways.
Let $v$ be a unitary in $M_n(D_{\phi}(A,B))$. Define $\phi^n:A \to \mathcal{L}_B(B^n)$ by $\phi^n(a)(b_1,b_2,\cdots,b_n)=(\phi(a)b_1,\phi(a)b_2,\cdots,\phi(a)b_n)$. Let $B^n\oplus B^n$ be graded by $(x,y)\mapsto (x,-y)$. Then
\[\left( B^n\oplus B^n,\left( \begin{matrix}
                           \phi^n & 0\\
                           0  & \phi^n \end{matrix}\right), \left( \begin{matrix}
                           0 & v\\
                           v^*  & 0 \end{matrix}\right)      \right) \]
 is a Kasparov $(A,B)$-module. The class of this module depends only on the class of $v$ in $K_1(D_{\phi}(A,B))$ so that the construction gives rise to a group homomorphism $\Omega:K_1(D_{\phi}(A,B)) \to \KK(A,B)$.
\begin{lem}\label{L:Thomsen}
Let $\phi$ an absorbing representation from $A$ to $\mathcal{L}(H_B)=M(B)$ where $B$ is a stable $C\sp*$-algebra. Then
$\Omega: K_1(D_{\phi}(A, B)) \to \KK(A,B)$ is an isomorphism.
\end{lem}
\begin{proof}
See \cite[Theorem 3.2]{Th}. In fact, Thomsen proved $K_1(\frac{D_{\phi}(A,B)}{(D_{\phi}(A,A,B))})$ is isomorphic to $\KK(A,B)$ via a map $\Theta$ where $D_{\phi}(A,A,B)=\{x\in D_{\phi}(A,B)\mid x\phi(A) \subset B\}$ is the ideal of $D_{\phi}(A,B)$. However, the same proof shows $\Omega$ is an isomorphism. Alternatively we can show that $K_i(D_{\phi}(A,A,B))=0$ for $i=0,1$ by the argument of \cite[Lemma 1.6]{H} with the fact that $K_*(M(B))=0$.
Thus, using the six term exact sequence,  $K_*(D_{\phi}(A,B))$ is isomorphic to $K_*\left(\frac{D_{\phi}(A,B)}{(D_{\phi}(A,A,B))}\right)$. This implies the map $\Omega$ which is the composition with $\Theta$ and $q_1$ is an isomorphism. Here $q_1$ is the induced map between K-groups from the quotient map from $D_{\phi}(A, B)$ onto $\frac{D_{\phi}(A,B)}{(D_{\phi}(A,A,B))}$.
\end{proof}
\begin{defn}\cite[Definition 3.2]{DE}
If $\pi,\sigma:A \to \mathcal{L}(E)$ are representations, we say that $\pi$ and $\sigma$ are properly asymptotically unitarily equivalent  and write $\pi \approxeq \sigma$ if  there is a continuous path of unitaries $u: [0,\infty ) \to \mathcal{U}(\mathcal{K}(E)+\mathbb{C}I_E)$, $u=(u_t)_{t\in [0,\infty )}$ such that for all $a\in A$
\begin{itemize}
\item[(i)] $lim_{t\to \infty}\| \sigma(a)-u_t \pi(a)u_t^*\|=0$,\\
\item[(ii)] $ \sigma(a)-u_t \pi(a)u_t^* \in \mathcal{K}(E)$ for all $t\in [0,\infty ).$
\end{itemize}
\end{defn}
In the above, we introduced the Fredholm picture of $\KK$-group. There is an alternative way to describe the element of $KK$-group. The Cuntz picture is described by a pair of representations $\phi,\psi:A \to \mathcal{L}(H_B)=M(B\otimes K)$ such that $\phi(a)-\psi(a)\in \mathcal{K}(H_B)=B\otimes K$. Such a pair is called a Cuntz pair. They form a set denoted by $\mathbb{E}_h(A,B)$. A homotopy of Cuntz pairs consists of a Cuntz pair $(\Phi,\Psi):A \to M(C([0,1])\otimes (B\otimes K))$. The quotient of $\mathbb{E}_h(A,B)$ by homotopy equivalence  is a group $\KK_h(A,B)$ which is isomorphic to $\KK(A,B)$ via the mapping sending $[\phi,\psi]$ to $[\phi,\psi,1]$.

 Dadarlat and Eilers proved that $[\phi,\psi]=0$ in $\KK_h(A,B)$ if and only if there is a representation $\gamma:A \to M(B\otimes K)=\mathcal{L}(H_B)$ such that $\phi\oplus\gamma \approxeq \psi\oplus\gamma$ \cite[Proposition 3.6]{DE}. The point is that the equivalence is implemented by unitaries of the form  compact $+$ identity. Sometimes, we can have a non-stable equivalence keeping this useful point.

\begin{defn}
 Let  $A$ be a $C\sp*$-algebra. Denote by ${\widetilde A}$ its unitization. We say that $A$ has $K_1$-injectivity if the map from $\mathcal{U}(\widetilde{A})/\mathcal{U}_0(\widetilde{A})$ to $K_1(A)$ is injective where $\mathcal{U}(\widetilde{A})$ is the unitary group and $\mathcal{U}_0(\widetilde{A})$ is the connected component of the identity. We note that H. Lin proved in \cite[Lemma 2.2]{Lin96} that real rank zero implies $K_1$-injectivity.
\end{defn}

\begin{thm}\label{T:DE}
Let $A$ be a separable $C\sp*$-algebra and let $\psi, \phi: A \to H_B$ be a Cuntz pair of absorbing representations. Suppose that the composition of $\phi$ with the natural quotient map $\pi:M(B\otimes K) \to M(B\otimes K)/ B\otimes K$, which will be denoted by $\dot{\phi}$, is faithful. Further, we suppose that $D_{\phi}(A,B)$ satisfies $K_1$-injectivity. If $[\phi,\psi]=0$ in $KK(A,B)$, then $\phi \approxeq \psi$.
\end{thm}
\begin{proof}
The proof of this theorem is almost identical to the one given in \cite[Theorem 3.12]{DE}. We just give the proof to illustrate how our assumptions play the roles.

By Theorem \ref{T:Voiculescu}, we get a continuous family of unitaries $(u_t)_{t\in\left[0,\infty\right)}$ in $M(B\otimes K)$ such that
\begin{equation}\label{E:asymt}
u_t\phi(a)u_t^*- \psi(a) \in C_0(\left[0,\infty \right))\otimes (B\otimes K).
\end{equation}
Note that (\ref{E:asymt}) implies $[\phi,\psi]=[\phi,u_1\phi u_1^*]$ (See \cite[Lemma 3.1]{DE}).
We assume that $[\phi,\psi]=0$ and we conclude that $[\phi,u_1 \phi u_1^*]=0$. Since $(\phi,\phi,u_1^*)$ is unitarily equivalent to $(\phi,u_1\phi u_1^*,1)$, \[[(\phi,\phi,u_1)]=[(\phi,\phi,u_1^*)]=0.\]
Since the isomorphism $\Omega:K_1(D_{\phi}(A,B)) \to \KK(A,B)$  sends $[u_1]$ to $[\phi,\phi,u_1]$  by Lemma \ref{L:Thomsen}, $K_1$-injectivity implies that $u_1$ is homotopic to $1$ in $D_{\phi}(A,B)$. Thus we may assume that $u_0=1$ in (\ref{E:asymt}).

Let $E_{\phi}$ be a $C\sp*$-algebra $\phi(A)+B\otimes K$.
We define $(\alpha_t)_{t\in \left[ 0, \infty \right)}$ in $\Aut_{0}(E_{\phi})$ by $\Ad (u_t)$. Note that $\alpha_0=\id $ and $(\alpha_t)$ is a uniform continuous family of automorphisms. Thus we apply Proposition 2.15 in \cite{DE} and get a continuous family $(v_t)_{\left[ 0, \infty\right)}$ of unitaries in $E_{\phi}$ such that
\begin{equation}\label{E:Aut}
 \lim_{t \to \infty}\| \alpha_{t}(x) -\Ad v_t(x)\|=0
\end{equation}
for any $x \in E_{\phi}$.

Combining (\ref{E:Aut}) with (\ref{E:asymt}), we obtain $(v_t)_{\left[ 0, \infty\right)}$ of unitaries in $E_{\phi}$ such that
\[  \lim_{t \to \infty}\| v_t\phi(a) v_t^*-\psi(a)\|=0 \]
for any $a \in A$. Since $\dot{\phi}$ is faithful, we can replace  $(v_t)_{\left[ 0, \infty\right)}$ by a family of unitaries in
$B\otimes K +\mathbb{C}1$ by the argument shown in the STEP 1 of the proof of Proposition 3.6 in \cite{DE}.
\end{proof}
Recall the definition of the esssential codimension of Brown, Douglas, and Fillmore defined by  two projections $p,q$ in $B(H)$  whose difference is compact as we have defined  in Introduction. Using KK-theory, or K-theory, we generalize this notion as follows, keeping the same notation.
\begin{defn}\label{D:BDF}
Given two projections $p,q \in M(B\otimes K)$ such that $p-q \in B\otimes K$, we consider representations $\phi,\psi $ from $\mathbb{C}$ to $M(B\otimes K)$ such that $\phi(1)=p,\psi(1)=q$. Then $(\phi,\psi)$ is a Cuntz pair so that
we define $[p:q]$ as the class $[\phi,\psi]\in \KK(\mathbb{C},B) \simeq K(B)$.
\end{defn}

\begin{lem}\cite{Lin1}\label{L:Lin}
Let $B$ be a non-unital ($\sigma$-unital) purely infinite simple $C\sp{*}$-algebra. Let $\phi, \psi$ be two monomorphisms from $C(X)$ to $M(B\otimes K)$ where $X$ is a compact metrizable space. If $\dot{\phi},\dot{\psi}$ are still injective, then they are approximately unitarily equivalent.
\end{lem}
The following theorem is a sort of generalization of BDF's result about the essential codimension.
\begin{thm}\label{T:BDF}
Let $B$ be a non-unital ($\sigma$-unital) purely infinite simple $C\sp{*}$-algebra such that $M(B\otimes K)$ has real rank zero.
Suppose two projections $p$ and $q$ in $M(B\otimes K)=\mathcal{L}(H_B)$ such that $p-q \in B\otimes K$ and neither of them are in $B\otimes K$. If $[p,q] \in K_{0}(B)$ vanishes, then there is a unitary $u$ in $ \id + B\otimes K$ such that $upu^*=q $.
\end{thm}
\begin{proof}
Step1: Let $\phi,\psi:\mathbb{C} \to M(B\otimes K)$ be representations from $p$ and $q$ respectively.
Evidently $\phi$ is injective. Moreover, it does not contain any ``compacts'' since $p$ does not belong to $B\otimes K$. Thus $\dot{\phi}$ is faithful.
Recall $\psi_{\infty}$ is defined by $\psi_{\infty}(a)=\sum S_i \psi(a) S_i^*$ where $\{S_i\}$ is a sequence of isometries in $M(B\otimes K)$ such that $S_iS_j^*=0$ for $i\neq j$. Suppose that $\psi_{\infty}(\lambda)=0$ for $\lambda \in \mathbb{C}$. Then $S_i^*\psi_{\infty}(\lambda) S_i=\psi(\lambda)=0$ or  $\lambda q=0$. Thus $\lambda=0$. Similarly, $\dot{\psi}_{\infty}$ is injective.  Then they are approximately unitarily equivalent by  applying Lemma \ref{L:Lin} to $X=\{x_0\}$. Thus we have a unitary $U$ in
$\mathcal{L}(H_B)$ such that
\begin{equation} \label{E:unitaryequivalecne}
U^* \phi(a)U - \psi_{\infty}(a)
\end{equation}
 for $a \in \mathbb{C}$.

Note that to get a sequnce of isometries $\{v_i\} \in \mathcal{L}(H_B^{(\infty)},H_B)$ satisfying the conditions of  Lemma \ref{L:absorbing}, what we needed was the equation (\ref{E:unitaryequivalecne}).  Following the same argument in the proof of  Theorem \ref{T:Voiculescu}, we get $\phi \underset{\text{asym}}{\sim} \psi$. In other words, we have a continuous family of unitaries $(u_t)_{t\in\left[0,\infty\right)}$ in $M(B\otimes K)$ such that
\begin{equation*}\label{E:asym}
u_t\phi(a)u_t^*- \psi(a) \in C_0(\left[0,\infty \right))\otimes (B\otimes K) \quad\text{for any $a$ in A}.
\end{equation*}
Since $D_{\phi}(\mathbb{C},B)$ has real rank zero, it satisfies $K_1$-injectivity. Thus it follows that $\phi \approxeq \psi$ as in the proof of Theorem \ref{T:DE}.

Step2: For large enough $t$, we can take $u_t=u$ of the form `identity + compact' such that $\|upu^*-q\| <1$. For the moment we write $upu^*$ as $p$. Thus $\|p-q\| <1$.
Note that $p-q \in B\otimes K$. Then  $z=pq+(1-p)(1-q)\in 1+ B\otimes K $ is invertible and
 $pz=zq$. If we consider the polar decomposition of $z$ as $z=v|z|$.
 It is easy to check that  $v \in 1+B\otimes K$ and $vpv^{\ast}=q$. Now $w=vu$ is also a unitary of the form `identity + compact' such that
 \[ wpw^*=q.\]

\end{proof}

\section{Application: projection lifting}
In this section, we show an application of proper asymptotic unitary equivalence of two projections. In this application, with an additional real rank zero property, the unitary of the form `identity + compact' plays a crucial role as we shall see.

Let $B$ be a stable $C\sp*$-algebra such that the multiplier algebra $M(B)$ has real rank zero.
Let $X$ be $[0,1], [0,\infty)$, $(-\infty,\infty)$ or
$\mathbb{T}=[0,1]/\{0,1\}$. When $X$ is compact, let $I=C(X)\otimes
   B$ which is the $C\sp{*}$-algebra of (norm continuous) functions from $X$ to $B$. When $X$ is not compact,  let $I=C_0(X)\otimes B$ which is the $C\sp*$-algebra of continuous functions from $X$ to $B$ vanishing at infinity. Then $M(I)$ is given by $C_b(X, M(B)_s)$, which is the set of
   bounded functions from $X$ to $B(H)$, where $M(B)$ is given the strict
   topology. Let $\mathcal{C}(I)=M(I)/I$ be the corona algebra of $I$ and also
   let $\pi:M(I) \to \mathcal{C}(I)$ be the natural quotient map. Then an element
   $\mathbf{f}$ of the corona algebra can be represented as follows:  Consider a
 finite partition of $X$, or $X \smallsetminus \{0,1\}$ when $X=\mathbb{T}$ given by partition points $x_1 < x_2 < \cdots
 < x_n $ all of which are in the interior of $X$ and divide $X$ into
 $n+1$ (closed) subintervals $X_0,X_1,\cdots,X_{n}$. We can take $f_i \in
 C_b(X_i, M(B)_s)$ such that $f_i(x_i) -f_{i-1}(x_i)\in B$
 for $i=1,2,\cdots,n$ and $f_0(x_0)-f_n(x_0) \in B$  where $x_0=0=1$ if $X$ is $\mathbb{T}$.

 \begin{lem}
 The coset in $\mathcal{C}(I)$ represented by
 $(f_0,\cdots,f_n)$  consists of functions $f$ in $M(I)$ such that $f- f_i \in
 C(X_i)\otimes B$ for every $i$ and $f-f_i $ vanishes (in norm) at any
 infinite end point of $X_i$.
 \end{lem}
\begin{proof}
 If $X$ is compact, then we set $x_0=0$, $x_{n+1}=1$. Otherwise, we set $x_0=x_1-1$ when $X$ contains $-\infty$, and $x_{n+1}=x_n +1$ when $X$ contaions $+\infty$.
Then we define a function in $C(X)\otimes B$ by
  \begin{equation*}
m_i(x)=
\begin{cases}
\frac{x-x_{i-1}}{x_i-x_{i-1}}(f_i(x_i) -f_{i-1}(x_i)), &\text{if $x_{i-1}\leq x \leq x_i$}\\
\frac{x-x_{i+1}}{x_i-x_{i+1}}(f_i(x_i) -f_{i-1}(x_i)), &\text{if $x_{i}\leq x \leq x_{i+1}$}\\
0,                               &\text{otherwise}
\end{cases}
\end{equation*}
for each $i=1,\cdots n$. In addition, we set $m_0=m_{n+1}=0$.
Then we define a function $\widetilde{f}$ from $f_i$'s by
\[\widetilde{f}(x)=f_i(x)-m_i(x)/2+m_{i+1}(x)/2\] on each $X_i$.
It follows that $f_i(x_i)-m_i(x_i)/2+m_{i+1}(x_i)/2=f_{i-1}(x_i)-m_{i-1}(x_i)/2+m_{i}(x_i)/2$. Thus $\widetilde{f}$ is well defined.
 The conditions $f-f_i \in C(X_i)\otimes B$ for each $i$ imply that $f-\widetilde{f}$ is  norm continuous function from $X$ to $B$ since
$f|_{X_i}(x_i) -\widetilde{f}|_{X_i}(x_i)= f|_{X_{i-1}}(x_i)-\widetilde{f}|_{X_{i-1}}(x_i)$.

\end{proof}
 Similarly $(f_0,\cdots,f_n)$ and $(g_0,\cdots,g_n)$
 define the same element of $\mathcal{C}(I)$ if and only if $f_i - g_i \in
 C(X_i)\otimes B$ for $i=0,\cdots,n$ if $X$ is compact.  $(f_0,\cdots,f_n)$ and $(g_0,\cdots,g_n)$
 define the same element of $\mathcal{C}(I)$ if and only if  $f_i - g_i \in
 C(X_i)\otimes B$ for $i=0,\cdots,n-1$, $f_n -g_n \in
 C_0([x_n,\infty))\otimes B$  if $X$ is $[0.\infty)$.  $(f_0,\cdots,f_n)$ and $(g_0,\cdots,g_n)$
 define the same element of $\mathcal{C}(I)$ if and only if
  $f_i - g_i \in
 C(X_i)\otimes B$ for $i=1,\cdots,n-1$, $f_n -g_n \in
 C_0([x_n,\infty))\otimes B$, $f_0-g_0 \in
 C_0((-\infty,x_1])\otimes B$ if $X=(-\infty,\infty)$.\\

The following theorem says that any projection in the corona algebra of $C(X)\otimes B$ for some $C\sp{*}$-algebras $B$ is described by a ``locally trivial fiber bundle'' with the fibre $H_{B}$ in the sense of Dixmier and Duady \cite{DixDua}.
\begin{thm}\label{T:locallift}
Let $I$ be $C(X)\otimes B$ or $C_0(X)\otimes B$ where $B$ is a stable $C\sp{*}$-algebra such that $M(B)$  has real rank zero.
Then a projection $\mathbf{f}$ in $M(I)/I$ can be represented by $(f_0,f_1,\cdots, f_n)$ as above where $f_i$ is a projection valued function in
$C(X_i)\otimes M(B)_s$ for each $i$.
\end{thm}
\begin{proof}
Let $f$ be the element of $M(I)$ such that $\pi(f)=\mathbf{f}$.
 Without loss of generality, we can assume $f$ is self-adjoint and $0 \leq f \leq 1$.
\begin{itemize}
 \item [(i)]Suppose $X$ does not contain any infinite point.
 Choose a point $t_0 \in X $. Then there is a self-adjoint element $T \in M(B)$ such that $T-f(t_0) \in B$ and
 the spectrum of $T$ has a gap around $1/2$ by \cite[Theorem 3.14]{BP}. So we consider $f(t)+T-f(t_0)$ which is still self-adjoint whose image is $\mathbf{f}$.
 Thus we may assume $f(t_0)$ is a self-adjoint element whose spectrum has a gap around $1/2$.

 Since $r(f(t)):t \to f(t)-f(t)^{2}$ is norm continuous where $r(x)=x-x^2$,
 if we pick a point $z$ in $\left( 0,\frac{1}{4}\right)$ such that $z \notin \sigma(f(t_0)-f(t_0)^2)$,
 then  $\sigma(f(s))$ omits $r^{-1}(J)$  for $s$ sufficiently
 close to $t$ where $J$ is an interval containing $z$. In other
 words, there is  $\delta>0$ and $b> a >0 $ such that if $|t_0-s|<\delta $, then $ \sigma(f(s)) \subset
 [0, a)\cup (b,1]$. \\
  If we let $f_{t_0}(s)=\chi_{(b,1]}(f(s))$ for $s$ in $(t_0 - \delta,
  t_0+\delta)$ where $\chi_{(b,1]}$ is the characteristic function on $(b,1]$, then it is a continuous projection valued function such that
  $f_{t_0} - f \in C(t_0-\delta,t_0+\delta)\otimes B$.

  By repeating the above procedure, since $X$ is compact, we can find $n+1$ points $t_0,\cdots, t_n$, $n+1$ functions $f_{t_0}, \cdots,  f_{t_n}$,
  and an open covering $\{O_i\}$ such that $t_i \in O_i$, $O_i \cap O_{i-1} \ne
  \varnothing$, and  $f_{t_i}$ is projection valued function on $O_i$. Now
   let $f_i=f_{t_i}$ as above. Take the point $x_i \in O_{i-1} \cap O_i$ for
  $i=1, \cdots, n$. Then $f_{i}(x_i) - f_{i-1}(x_i)= f_{i}(x_i) - f(x_i) +
  f(x_i)- f_{i-1}(x_i) \in B$ and $f_0(x_0)-f_n(x_0) \in B$ if applicable. Let $X_{i}=[x_{i},x_{i+1}]$ for $i=1,
  \cdots, n-1$, $X_0=[0,x_1]$, and $X_n=[x_n,1]$. Since each $f_i$ is also
  defined on $X_{i}$, $(f_0, \cdots, f_n)$ is what we
  want.
  \item[(ii)] let $X$ be $[0, \infty)$.  Since $f^2(t) - f(t) \to 0 $ as $t$
  goes to $\infty$, for given $ \delta$ in $(0, 1/2) $, there is $M >0$ such that
  whenever $t \geq M$ then $\|f^2(t)-f(t) \| < \delta -\delta^2 $. It follows
  that $\sigma(f(t)) \subset [0, \delta) \cup (1-\delta,1]$ for $t \geq M$. Then
  again $\chi_{(1-\delta,1]}(f(t))$ is a continuous projection valued function for $t \geq M$
  such that $f(t) - \chi_{(1-\delta,1]}(f(t))$ vanishes in norm as $t$ goes to $\infty$.  By applying
  the argument in $(i)$ to $[0,M]$, we get a closed sub-intervals $X_i$ for $i=0,
  \cdots, n-1$ of $[0,M]$ and $f_{i} \in C_{b}(X_{i}, B(H))$. Now if we let
  $X_n=[M, \infty)$ and $f_n(t) = \chi_{(1-\delta,1]}(f(t)) $, we are done.
  \item[(iii)]The case $X=(-\infty, \infty)$ is similar to (ii).
  \end{itemize}
\end{proof}

When  a projection $\mathbf{f}\in C(I)$ is represented by $(f_0,f_1,\cdots,f_n)$  by Theorem \ref{T:locallift}, we note that  $f_i(x)$ is a projection in $M(B\otimes K)$ for each $x \in X_i$ and $f_i(x_i)-f_{i-1}(x_i)\in B$. Applying Definition \ref{D:BDF} we have $K$-theoretical terms $k_i=[f_i(x_i):f_{i-1}(x_i)]\in \KK(\mathbb{C},B)$ for $i=1,2,\cdots,n$. The following theorem shows that if all $k_i$'s are vanishing, then  a projection $\mathbf{f}$ in $C(I)$ lifts to a projection in $M(I)$.
\begin{thm}\label{T:lifting}
Let $I$ be $C(X)\otimes B$ where $B$ is a $\sigma$-unital, non-unital, purely infinite simple $C\sp{*}$-algebra such that $M(B)$  has real rank zero or $K_1(B)=0$ (See \cite{Zh}).
Let a projection $\mathbf{f}$ in $M(I)/I$ be represented by $(f_1,f_2,\cdots, f_n)$, where $f_i$ is a projection valued function in
$C(X_i)\otimes M(B)_s$ for each $i$, as in Theorem \ref{T:locallift}.
If $k_i=[f_i(x_i):f_{i-1}(x_i)]=0$ for all $i$, then the projection $\mathbf{f}$ in $M(I)/I$ lifts.
\end{thm}
\begin{proof}
Note that, by Zhang's dichotomy, $B$ is stable \cite[Theorem 1.2]{Zh}.
By induction, assume that $f_{j}(x_j)=f_{j-1}(x_j)$ for $j=1,2, \cdots, i-1$.

Let $f_i(x_i)=p_i, f_{i-1}(x_i)=p_{i-1}$. Since $[p_i:p_{i-1}]=0$, we have a unitary $u$ of the form ` identity + compact ' such that
$\|p_{i}-u^*p_{i-1}u\|<1/2$ by Theorem \ref{T:BDF}. Since $B$ has real rank zero,  given $0< \epsilon <1/4 $ there is a unitary $v \in \mathbb{C}1+B$ with finite spectrum such that $\|u-v \|<\epsilon$ \cite{Lin93}, \cite{Lin96}. Then \[ \|p_{i}-vp_{i-1}v^* \| \leq \|p_{i}-up_{i-1}u^* \|+\|up_{i-1}u^*-vp_{i-1}v^*\|<1.\]
Note that  $p_{i}-vp_{i-1}v^* \in B$. Thus we have $wp_{i}w^*=vp_{i-1}v^*$ for some unitary $w \in \id+B$. (Recall that Step 2 of the proof of Theorem \ref{T:BDF}.) Let $g_i=wf_iw^*$, then $f_i-g_i\in C(X_i)\otimes B $ since $w$ is of the form ` identity + compact'.

On the other hand, we can write $v$ as $e^{ih}$ where $h$ is a self-adjoint element in $B$ since $v$ has the finite spectrum. A homotopy of unitaries $t \to e^{ith}$, which are of the form `` identity + compact'', connects $1$ to $v$. Now we define $g_{i-1}(t)$ as
 \[\exp \left(i \frac{t-x_{i-1}}{x_i-x_{i-1}}h\right)f_{i-1}(t)\overline{\exp\left(i\frac{t-x_{i-1}}{x_i-x_{i-1}}h\right)}\] for $t \in [x_{i-1},x_i]$. Then
we see that $g_{i-1}(x_i)=g_{i}(x_i)$, $g_{i-1}-f_{i-1}\in C(X_{i-1})\otimes K$, and $g_{i-1}(x_{i-1})=f_{i-1}(x_{i-1})$. Moreover, if we let $g_{i+1}=wf_{i+1}w^*$, then $f_{i+1}-g_{i+1} \in C(X_{i+1})\otimes B $, and
\[
\begin{split}
[g_{i+1}(x_{i+1}):g_{i}(x_{i+1})]&=[wf_{i+1}(x_{i+1})w^*:wf_i(x_i+1)w^*]\\
&=[f_{i+1}(x_{i+1}):f_i(x_{i+1})]=0.
\end{split}
\]
Then $(f_0,f_1,\cdots, f_n)$ and  $(f_0,f_1,\cdots,g_{i-1},g_i,g_{i+1},f_{i+2},\cdots, f_n)$  define the same element  $\mathbf{f}$ while the $k_i$'s are unchanged and $i$-th discontinuity is resolved. So we take the latter as $(f_0,\cdots,f_n)$ such that $f_j(x_j)=f_{j-1}(x_j)$ for $j=1,\dots,i$.
We can repeat the same procedure until we have $f_i(x_i)=f_{i-1}(x_i)$ for all $i$. It follows that $(f_0,\cdots,f_n)$ is a projection in $M(C(X)\otimes B)$ which lifts $\mathbf{f}$.
\end{proof}
\begin{rem}
When $I=C_0(X)\otimes B$ where $X$ is $ [0,\infty)$ or $(-\infty,\infty)$, the similar result holds replacing $C(X_i)\otimes B$ with $C_0(-\infty,x_1]\otimes B$ or $C_0[x_n,\infty) \otimes B$  for $i=0$ or $i=n$ respectively.
\end{rem}
\section{Acknowledgements}
 Although this work was not carried out at Purdue, a significant influence on the auther was made by Larry Brown and Marius Dadarlat  who have aquainted him with geometric ideas in operator algebras.
He also would like to thank Huaxin Lin for answering the question related to Lemma 2.14.

\end{document}